# Novel Modifications of Parallel Jacobi Algorithms[*]

Sanja Singer[†], Saša Singer[‡], Vedran Novaković[§],
Aleksandar Ušćumlić,[¶]
and Vedran Dunjko[‖]

September 26, 2018

## Abstract

We describe two main classes of one-sided trigonometric and hyperbolic Jacobi-type algorithms for computing eigenvalues and eigenvectors of Hermitian matrices. These types of algorithms exhibit significant advantages over many other eigenvalue algorithms. If the matrices permit, both types of algorithms compute the eigenvalues and eigenvectors with high relative accuracy.

We present novel parallelization techniques for both trigonometric and hyperbolic classes of algorithms, as well as some new ideas on how pivoting in each cycle of the algorithm can improve the speed of the parallel one-sided algorithms. These parallelization approaches are applicable to both distributed-memory and shared-memory machines.

The numerical testing performed indicates that the hyperbolic algorithms may be superior to the trigonometric ones, although, in theory, the latter seem more natural.

## 1 Introduction

Among a variety of diagonalization algorithms for Hermitian and symmetric matrices, the Jacobi algorithm is the oldest and the simplest one, but is often considered too slow for practical usage. However, Jacobi-type algorithms have recently returned to the focus of active research, mostly due to their accuracy properties and inherent amenability for parallelization.

Demmel and Veselić [8] showed that both the one-sided and the two-sided Jacobi algorithms are accurate in the relative sense for positive definite matrices. To explain precisely what the "accuracy in relative sense" means, we need to define the scaled condition $\kappa_s(H)$ of a Hermitian positive definite matrix $H$. Let $H$ be scaled such that

$$A := D^{-1} H D^{-1}, \qquad D = \operatorname{diag}((h_{11})^{1/2}, \ldots, (h_{nn})^{1/2}), \qquad (1.1)$$

[*]This work was supported by grant 037–1193086–2771 by the Ministry of Science, Education and Sports, Republic of Croatia.

[†]Faculty of Mechanical Engineering and Naval Architecture, University of Zagreb, Ivana Lučića 5, 10000 Zagreb, Croatia, e-mail: `ssinger@fsb.hr`

[‡]Department of Mathematics, University of Zagreb, P.O. Box 335, 10002 Zagreb, Croatia, e-mail: `singer@math.hr`

[§]Faculty of Mechanical Engineering and Naval Architecture, University of Zagreb, Ivana Lučića 5, 10000 Zagreb, Croatia, e-mail: `venovako@fsb.hr`

[¶]MSV sustavi d.o.o., Tatjane Marinić 12, 10430 Samobor, Croatia, e-mail: `sasa.uscumlic@gmail.com`

[‖]School of EPS – Physics Department, David Brewster Building, Heriot-Watt University, Edinburgh, EH14 4AS, United Kingdom, e-mail: `vd51@hw.ac.uk`

where $n$ is the order of $H$. Then $\kappa_s(H) := \kappa_2(A) = \|A\|_2 \|A^{-1}\|_2$. According to [23], $A$ is nearly optimally scaled, in the sense that

$$\kappa_2(A) \leq n \min_\Delta \kappa_2(\Delta H \Delta),$$

over all diagonal matrices $\Delta$. The previous inequality implies $\kappa_2(A) \leq n\kappa_2(H)$, but it frequently happens that $\kappa_2(A) \ll \kappa_2(H)$.

Demmel and Veselić in [8] also proved that the eigenvalues $\lambda_i$ of $H$ can be computed with a small relative error, i.e., they replaced the standard relative bound

$$\frac{|\lambda_i' - \lambda_i|}{\lambda_i} \leq \varepsilon f(n) \kappa_2(H),$$

by

$$\frac{|\lambda_i' - \lambda_i|}{\lambda_i} \leq \varepsilon f(n) \kappa_2(A), \tag{1.2}$$

where $\varepsilon$ is the machine precision, and $f$ is a slowly increasing function of the order $n$.

In the case of a positive definite $H$, the two-sided Jacobi algorithm can be viewed as the singular value decomposition (SVD) of a factor $G$ of the matrix $H$

$$H = GG^*.$$

Moreover, for a positive definite $H$, the eigenvalues of $G^*G$ and $GG^*$ are equal, so we can choose to transform either $GG^*$ or $G^*G$. The sequence of two-sided unitary transformations which diagonalizes $GG^*$ needs to be applied only from one side, say the right-hand side, i.e., on $G^*$.

All the angles in this process are calculated from the temporary iterate $H_\ell := G_\ell G_\ell^*$, but applied only on $G_\ell^*$. If the iterates $G_\ell G_\ell^*$ tend to a diagonal matrix for $\ell \to \infty$, then $G_\ell^*$ tends to a matrix with orthogonal (but not orthonormal!) columns. The squared norms of the columns of the final $G_\ell^*$ are the eigenvalues of $H$. A similar fact holds for the matrix $G^*G$ and its factor $G$.

If $H \in \mathbb{C}^{n \times n}$ is an indefinite matrix of rank $m$, $m \leq n$, the diagonalization task is harder to deal with. The Jacobi-type algorithms that possess the relative accuracy property work on a factor of $H$. The factor of $H$ is computed by using Slapničar's modification (see [21]) of the Bunch–Parlett Hermitian indefinite factorization with pivoting ([1, 2, 3, 4, 5, 6])

$$\widehat{P} H \widehat{P}^T = GJG^*, \quad J = \text{diag}(j_{11}, \ldots, j_{mm}), \tag{1.3}$$

where $\widehat{P}$ is a permutation matrix, $G \in \mathbb{C}^{n \times m}$ is a block lower trapezoidal matrix with diagonal blocks of order one or two, and $j_{ii} \in \{-1, 1\}$, for $1 \leq i \leq m$. If $H$ is nonsingular, $G$ is a block lower triangular matrix.

Similarly to the positive definite case, the indefinite Jacobi diagonalization can be viewed as the hyperbolic SVD of the matrix $G$ from (1.3) with respect to the signature matrix $J$ (see, for example, [26]),

$$G = U\Sigma V^*,$$

where $U$ is an orthogonal matrix, $\Sigma$ is a diagonal matrix with non-negative entries and $V^*$ is a $J$-orthogonal matrix, i.e., $V^* J V = J$.

This SVD can be obtained either by orthogonal transformations applied to $G$ from the left, or by hyperbolic ones applied from the right. The former case is known as the one-sided Jacobi algorithm [25, 9], and the latter is known as the hyperbolic one-sided Jacobi algorithm [24].

If the factor $G$ is well-scaled, then both algorithms are accurate in the relative sense [22, 9]. Slapničar [22] generalized the proof of the relative accuracy [8] to the case of the hyperbolic Jacobi algorithm. Namely, if $H$ is a nonsingular indefinite matrix, and the relation (1.1) is replaced by

$$A = D^{-1} \widehat{H} D^{-1}, \qquad D = \text{diag}((\hat{h}_{11})^{1/2}, \ldots, (\hat{h}_{nn})^{1/2}),$$

where $\widehat{H} = \sqrt{H^2}$ is the positive definite polar factor of $H$, then the accuracy of the hyperbolic Jacobi algorithm is essentially given by (1.2). From the bounds for the one-sided Jacobi algorithms, it is obvious that the matrix $H$ permits an accurate computation of eigenvalues if the scaled condition $\kappa_2(A)$ is small.

On the other hand, Dopico, Koev and Molera [9] proved that the one-sided trigonometric Jacobi algorithm computes the eigenvalues of $H := X\widehat{D}X^*$, where $\widehat{D}$ is a diagonal nonsingular matrix, with a relative bound given by

$$\frac{|\lambda'_i - \lambda_i|}{\lambda_i} \leq O(\varepsilon \kappa_2(X)).$$

In this paper we develop a sequence of modifications, applicable to both trigonometric and hyperbolic one-sided algorithms, aimed at increasing the speed of parallel implementations of these algorithms. The combined effect is a speedup of approximately 30% over the straight-forward parallel realizations of the trigonometric and hyperbolic [20] algorithms.

The rest of the paper is organized as follows. In Section 2 we briefly describe the sequential Jacobi algorithms with emphasis on the details needed in a parallel implementation. Section 3 is devoted to detailed descriptions of the corresponding parallel algorithms. Specially, we discuss the advantages and the drawbacks of the algorithms, and present the modifications mentioned above. In the final section, we give some of the results of the performed numerical testing. In the Appendix we derive the error bounds for the eigenvector computation in the trigonometric case.

## 2 Block algorithms, pivot strategies and parallelization

### 2.1 Pointwise algorithms

The pointwise Jacobi algorithms, that diagonalize a $2 \times 2$ matrix in each step, both the two-sided and the one-sided ones, are well known and described in details, for instance, in [24, 22, 9, 13, 20].

The one-sided trigonometric algorithm operates on $G^*$ from the right, choosing pivots from $H = GJG^*$ (see [9]). The hyperbolic algorithm operates from the right on $G$, diagonalizing the matrix pair $(A, J) := (G^*G, J)$ (see [24]).

To unify the notation for both the trigonometric and the hyperbolic Jacobi-type algorithms, let $A^\circ = H$, $G^\circ = G^*$ in the trigonometric, and $A^\circ = A$, $G^\circ = G$ in the hyperbolic algorithm. Fig. 1 shows the relation of the matrix $A^\circ$ to its factor $G^\circ$. The matrix $A^\circ_P$ is called the pivot (sub)matrix.

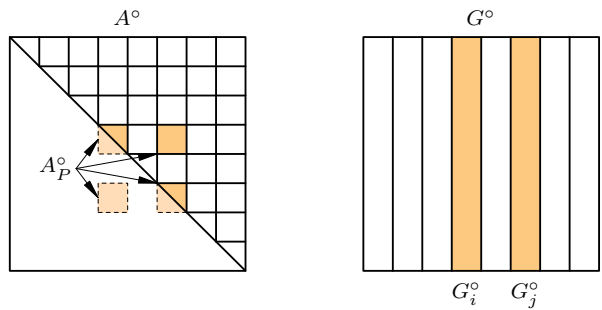

Figure 1: Block columns $G^\circ_P := [G^\circ_i, G^\circ_j]$ in $G^\circ$ correspond to a square block $A^\circ_P$ in $A^\circ$.

We summarize the pointwise trigonometric and hyperbolic algorithms below. Both algorithms operate on a chosen pair $G^\circ_P := [g^\circ_i, g^\circ_j]$ of ordinary columns of $G^\circ$.

Trigonometric Jacobi:

1. Find the next pivot pair $(i,j)$ and compute
$$A_P^\circ := H_P = \begin{bmatrix} h_{ii} & h_{ij} \\ h_{ij}^* & h_{jj} \end{bmatrix} = \begin{bmatrix} g_i \\ g_j \end{bmatrix} J \begin{bmatrix} g_i^* & g_j^* \end{bmatrix}.$$

2. Diagonalize $A_P^\circ$ by a single trigonometric rotation,
$$Q_P = \begin{bmatrix} \cos\varphi & e^{i\alpha}\sin\varphi \\ -e^{-i\alpha}\sin\varphi & \cos\varphi \end{bmatrix}, \quad (2.1)$$
i.e., find $\cos\varphi$, $\sin\varphi$ and $e^{i\alpha}$ such that
$$Q_P^* A_P^\circ Q_P = \begin{bmatrix} h_{ii}' & 0 \\ 0 & h_{jj}' \end{bmatrix}.$$

3. Apply the rotation to the columns $g_i^*$ and $g_j^*$,
$$\begin{bmatrix} (g_i^*)' & (g_j^*)' \end{bmatrix} = \begin{bmatrix} g_i^* & g_j^* \end{bmatrix} Q_P.$$

Hyperbolic Jacobi:

1. Find the next pivot pair $(i,j)$ and compute
$$A_P^\circ := A_P = \begin{bmatrix} a_{ii} & a_{ij} \\ a_{ij}^* & a_{jj} \end{bmatrix} = \begin{bmatrix} g_i^* \\ g_j^* \end{bmatrix} \begin{bmatrix} g_i & g_j \end{bmatrix}.$$

2. If the signs in $J_P$ (the $2\times 2$ submatrix of $J$ that corresponds to $A_P^\circ$) are the same, choose the trigonometric rotation (2.1), otherwise choose the hyperbolic rotation
$$Q_P = \begin{bmatrix} \cosh\varphi & e^{i\alpha}\sinh\varphi \\ e^{-i\alpha}\sinh\varphi & \cosh\varphi \end{bmatrix},$$
to diagonalize $A_P^\circ$,
$$Q_P^* A_P^\circ Q_P = \begin{bmatrix} a_{ii}' & 0 \\ 0 & a_{jj}' \end{bmatrix}.$$

3. Apply the rotation to the columns $g_i$ and $g_j$,
$$\begin{bmatrix} g_i' & g_j' \end{bmatrix} = \begin{bmatrix} g_i & g_j \end{bmatrix} Q_P.$$

## 2.2 Blocked algorithms

In order to speed up the algorithms, a pivot pair of columns should be replaced by a pivot block that contains a pair of block columns (see Fig. 1). The generalization of the pointwise algorithm to a blocked one transforms the whole pivot block in each step, and this can be done in two different ways:
- by reducing its off-diagonal norm (the so-called *block-oriented* algorithm) [13];
- by diagonalizing it (the so-called *full block* algorithm) [14].

The main steps of the block-oriented and the full block algorithm are as follows.

Block-oriented algorithm:

1. Before each block sweep, for all block columns $G_i^\circ$ form the matrix $A_D^\circ$,
$$A_D^\circ := \begin{cases} G_i J G_i^*, & \text{trigonometric case,} \\ G_i^* G_i, & \text{hyperbolic case.} \end{cases}$$
Annihilate the upper triangle of $A_D^\circ$ only once. Apply the accumulated rotations to $G_i^\circ$.

2. Find the next pair of pivot block columns $G_i^\circ$ and $G_j^\circ$. Form the pivot block $A_P^\circ$,
$$A_P^\circ = \begin{bmatrix} A_{ii}^\circ & A_{ij}^\circ \\ (A_{ij}^\circ)^* & A_{jj}^\circ \end{bmatrix},$$
and annihilate each element of $A_{ij}^\circ$ only once.

3. Apply the accumulated rotations to the block columns $G_i^\circ$ and $G_j^\circ$.

Full block algorithm:

1. Find the next pair of pivot block columns $G_i^\circ$ and $G_j^\circ$. Form the pivot block $A_P^\circ$,
$$A_P^\circ = \begin{bmatrix} A_{ii}^\circ & A_{ij}^\circ \\ (A_{ij}^\circ)^* & A_{jj}^\circ \end{bmatrix},$$
and diagonalize it.

2. Apply the accumulated rotations to the block columns $G_i^\circ$ and $G_j^\circ$.

The pivot matrices $A_D^\circ$ and $A_P^\circ$ are processed in a one-sided manner, i.e., they are first factorized, and then the corresponding pointwise one-sided algorithm is applied to the factor.

## 2.3 Block factorizations

A naïve implementation of the blocked Jacobi algorithm would apply each transformation directly to the tall and skinny block

$$G_P^\circ := \begin{bmatrix} G_i^\circ & G_j^\circ \end{bmatrix}$$

(a factor of $A_P^\circ$), which is very inefficient. In practice, each block $G_P^\circ$ is preprocessed to obtain a square factor $R_P^\circ$ of $A_P^\circ$. Then, $R_P^\circ$ is transformed, and these transformations are accumulated in a "work" matrix $Q_P^\circ$. Finally, $G_P^\circ$ is postmultiplied (only once) by the accumulated $Q_P^\circ$, so the "update" is now a BLAS 3 operation.

This kind of an accumulated application of transformations influences the overall algorithm speed tremendously. In principle, there are two ways to compute the required square factor $R_P^\circ$.

1. By forming $A_P^\circ$ and computing a suitable Hermitian factorization afterwards, i.e., the Cholesky factorization in the hyperbolic case

$$A_P^\circ = R^* R,$$

where $R$ is an upper triangular matrix, and the Hermitian indefinite factorization in the trigonometric case

$$A_P^\circ = P^T R^* J_P R P, \tag{2.2}$$

where $P$ is a permutation matrix, $R$ is a block upper triangular matrix with diagonal blocks of order 1 or 2, and $J_P$ is a diagonal signature matrix containing the inertia of $A_P^\circ$. Handling the extra permutation in this case is detailed in Section 3.

2. By the QR-like factorization of $G_P^\circ$, i.e., the ordinary QR factorization in the hyperbolic, and the hyperbolic QR factorization [18] in the trigonometric algorithm.

The former case is faster, as it involves the multiplication of two block columns and the factorization of a relatively small square matrix, while the QR-like approach has a significant overhead of preserving the input block by copying, and applying the computed transformations to the tall and skinny original matrix. The details of preprocessing of $A_P^\circ$ blocks will be described in Section 3, for each type of the Jacobi algorithms.

## 2.4 The parallelization

The two different pivot blocks,

$$G_P^\circ := \begin{bmatrix} G_i^\circ & G_j^\circ \end{bmatrix} \quad \text{and} \quad G_Q^\circ := \begin{bmatrix} G_k^\circ & G_l^\circ \end{bmatrix},$$

with $i \neq j$ and $k \neq l$, can be independently and simultaneously transformed if

$$\{i, j\} \cap \{k, l\} = \emptyset.$$

This property, together with blocking, is the basis of parallel implementations of these algorithms, in a sense that the independent blocks are assigned as independent tasks to computational processes.

The iterations of the Jacobi algorithm are usually called *sweeps* (or cycles). In the block algorithms, we distinguish block (or outer) sweeps over pivot blocks, and inner sweeps of the pointwise algorithm inside each of the pivot blocks.

A *block (outer) pivot strategy* is the order in which the pivot blocks $A_P^\circ$ are chosen. An *inner pivot strategy* is the order of annihilations inside all pivot blocks.

In a block sweep, the pivot blocks are processed according to the block-oriented or the full block algorithm. Recall, in an inner sweep of the block-oriented algorithm, each off-diagonal element of $A^\circ$ is annihilated only once (or finitely many times, if the chosen pivot strategy is quasicyclic, see [10, 12, 17]). In the full block algorithm each pivot block $A_P^\circ$ is diagonalized.

There are many choices of pivot strategies, but only a few of them are suitable for parallel computation and provably convergent. Our choice of the outer strategy is the modulus pivot strategy, described in [16]. This strategy, with the row-cyclic inner strategy, is weakly equivalent to the block row-cyclic strategy (with the same inner strategy), and, therefore, ensures the convergence of the Jacobi algorithm, as shown in recent works [13, 14].

The pointwise modulus strategy simultaneously annihilates the elements of $A^\circ$ on each antidiagonal. If $n$ is even, and the antidiagonal is denoted by an even number, e.g., by 2 in Fig. 2, the modulus strategy annihilates only $n/2 - 1$ elements. If the additional element is also annihilated (presented in a lighter hue), the strategy is called the modified modulus strategy. The same principle is used for the corresponding block strategies that operate on blocks, instead of elements.

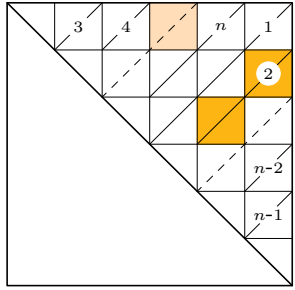

Figure 2: Annihilations of $A^\circ$ using the (modified) modulus pivot strategy.

Fig. 3 shows the block layouts for odd and even sweeps of the modulus pivot strategy. While the modulus pivot strategy can be inefficient in the sequential implementation (frequent cache spilling), it is ideal for parallel implementations.

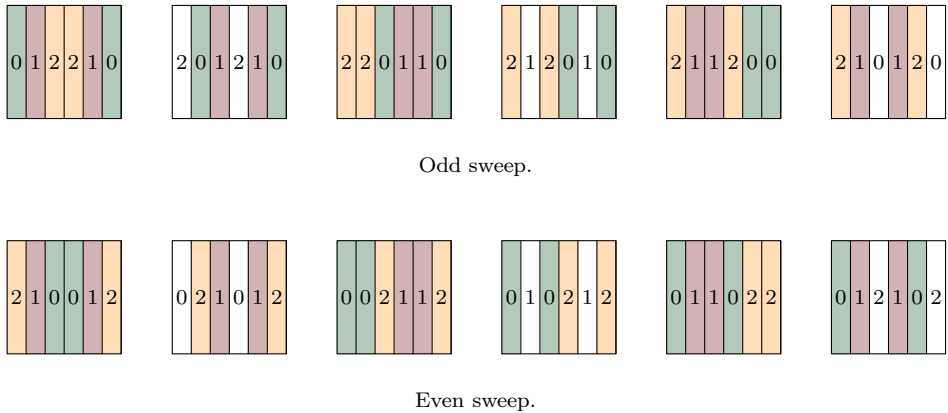

Figure 3: Modulus strategy for 6 blocks. A pair of block columns denoted by the same number label comprises a pivot block in each step. The white blocks are *skipped* (not transformed) in the modulus, but are *processed* in the modified modulus strategy.

# 3 Parallel implementations of the Jacobi algorithms

The trigonometric and the hyperbolic, both block-oriented and full block, Jacobi algorithms are parallelized in terms of single-threaded processes, communicating through the MPI (Message Passing Interface) stack. The implementations are independent of the underlying memory and network architecture, scaling from the symmetric multiprocessing (SMP), over the non-uniform memory architecture (NUMA), to the clusters of interconnected machines, provided only that the chosen MPI distribution efficiently supports the target architecture(s).

By $p \geq 2$ we denote the number of parallel processes involved in a running instance of any of our algorithms. These processes are arranged in a one-dimensional torus (ring), realized in MPI as a one-dimensional cyclic Cartesian virtual topology.

Initially, a block (with two block columns)

$$G_P^\circ := \begin{bmatrix} G_{r+1}^\circ & G_{2p-r}^\circ \end{bmatrix}$$

is assigned to a process with MPI rank $r$, where $0 \leq r < p$. The partitioning of $G^\circ$ into pivot blocks $G_P^\circ$ should be maintained as uniform as possible. In every step and at each process, the widths of the block columns in a block are allowed to differ by at most one, i.e., the widths of blocks $G_P^\circ$ of any two processes differ by at most two.

The ring topology is natural for the modulus pivot strategy (and vice versa), because it implies shifting of only one block column per process, cyclically after each step in a sweep. The cycle direction is defined with respect to the parity of a sweep number: in odd sweeps, the process $r$ sends a block column to the process $(r-1) \bmod p$, and receives one from the process $(r+1) \bmod p$, while in even sweeps, the process $r$ sends a block column to the process $(r+1) \bmod p$, and receives one from the process $(r-1) \bmod p$.

---

**Algorithm 3.1:** Modulus pivot strategy

**Initialization**: Process $r$ computes the indices of the initial column blocks $(i\_blk, j\_blk)$. An auxilliary pair of indices $(ip, jp)$ is used for the determination of the pivot indices in all subsequent steps:

$$ip = r + 1; \quad i\_blk = ip; \quad jp = 2 \cdot p - r; \quad j\_blk = jp.$$

**Description**: Routine `Next_Pair` computes the indices of the next pivot pair $(i\_blk, j\_blk)$ in the process $r$. It also tells whether to send/receive the block column $i\_blk$ or $j\_blk$. Routine `Send_Receive` determines the rank of the process from which the next block will be received, and the rank of the process to whom the next block will be sent.

`Next_Pair` $(r)$;
**begin**
    **if** $(ip + jp) > 2 \cdot p$ **then**
        $snd\_blk = i\_blk; \quad ip = ip + 1;$
        **if** $ip = jp$ **then**
            $ip = ip - p; \quad jp = ip;$
        **end if**
        $i\_blk = ip; \quad rcv\_blk = i\_blk;$
    **else**
        $snd\_blk = j\_blk; \quad jp = jp + 1;$
        $j\_blk = jp; \quad rcv\_blk = j\_blk;$
    **end if**
**end**

`Send_Receive` $(r)$;
**begin**
    **if** $(nsweep \bmod 2) > 0$ **then**
        $snd\_rnk = (p + r - 1) \bmod p;$
        $rcv\_rnk = (p + r + 1) \bmod p;$
    **else**
        $snd\_rnk = (p + r + 1) \bmod p;$
        $rcv\_rnk = (p + r - 1) \bmod p;$
    **end if**
**end**

---

To further elaborate the communication pattern in the modulus pivot strategy, let us describe the "process map" for a parallel Jacobi step. Instead of static processes and block columns being swapped among them, consider the block-partition of the matrix $A^\circ$, and assign a process to each block $A_P^\circ$ to be transformed in a given step of the (modified) modulus strategy. At the start of a sweep, all the processes lie on the antidiagonal blocks. When transitioning from one step to another, the processes move downwards, i.e., the row indices of the assigned blocks are incremented, while the column indices remain the same.

In an even step, a single (but each time different) process hits a diagonal block and changes for itself the above map traversing rule. The process keeps the column index of its block, and takes the row index unoccupied by other processes. In the subsequent steps, such a process keeps its block row index (say, $s$) fixed, and traverses the matrix $A^\circ$ row-wise, starting from the block $(s, s+1)$, and incrementing the block column index in each step. After the last step of a sweep the processes are positioned again on the antidiagonal, but their order from the start of a sweep is now reversed (see routine *Next_Pair* in Algorithm 3.1).

By looking at the process map we can tell whether a process replaces its first or its second block column in the after-step communication. While the process keeps moving downwards, it exchanges its first block column. When it hits a diagonal block, and afterwards, it exchanges its second block column. The efficient implementation of this traversing pattern and communication rules is given in the step transitioning procedure *Send_Receive* in Algorithm 3.1. The block column exchanges are realized with MPI, by using the MPI_SENDRECV routine (one could also choose MPI_SENDRECV_REPLACE).

## 3.1 Trigonometric Jacobi algorithms

We can now switch our focus to the detailed description of the operations performed by a single process, before the block column interchanges.

Formation of the pivot submatrix in the block case is the same as in the pointwise case, except that the ordinary columns $g_i^*$ and $g_j^*$ are now replaced by the block columns $G_i^*$ and $G_j^*$, owned by a particular process,

$$A_P^\circ := H_P = \begin{bmatrix} H_{ii} & H_{ij} \\ H_{ij}^* & H_{jj} \end{bmatrix} = \begin{bmatrix} G_i \\ G_j \end{bmatrix} J \begin{bmatrix} G_i^* & G_j^* \end{bmatrix}. \qquad (3.1)$$

The nonsingular pivot submatrix $H_P$ is then factored by the Bunch-Parlett Hermitian indefinite factorization with complete pivoting [6, 21], as in (2.2).

Note that the pivoting can change the affiliation of an individual column from one block column to the other. To prevent this, and to ensure the convergence of the algorithm, the columns should be restored to their original positions after the factorization.

If we apply the sequential one-sided trigonometric Jacobi algorithm to the factor $F = RP$ of the original $H_P$, i.e.,

$$H_P = F^* J_P F = (P^T R^*) J_P (RP),$$

we get the trigonometric full block (TF) and the trigonometric block-oriented (TB) algorithms. The similar naming convention will be used for all algorithms: the first letter denotes the type of the algorithm (trigonometric/hyperbolic), the second one denotes the type of blocking (block-oriented/full block), while the remaining letters describe the pivoting strategy.

If the one-sided Jacobi algorithm is applied directly to $R$, instead of $F$, the original column arrangement is lost, but some speedup is gained, especially in the full block (TFC) case. After the transformation of $R$ (either one sweep, or the full diagonalization) the rows of the unitary matrix $U_P$, that transforms $R$, are reordered according to the permutation used in the Bunch-Parlett factorization.

Though the convergence of TFC and TBC algorithms is not yet proven, we firmly believe that, due to the nature of complete pivoting in the Bunch-Parlett factorization, no column exchanges are needed after finitely many sweeps. As a result of pivoting, after the completion of the algorithm, the eigenvalues are sorted decreasingly in their absolute values.

However, there are two possible drawbacks in blocked trigonometric algorithms. The first one is that, even if $H$ is nonsingular, the pivot submatrix $H_P$ need not be.

**Example 1.** *Let*

$$H = \begin{bmatrix} 1 & 0 & 1 & 1 \\ 0 & 1 & 1 & 1 \\ 1 & 1 & 1 & 0 \\ 1 & 1 & 0 & 1 \end{bmatrix}.$$

The matrix $H$ is of rank 4, with eigenvalues $-1$, $1$, $1$, $3$. The Hermitian indefinite factorization of $H$ (with complete pivoting) gives $H = P^T R^* J R P$, where

$$R = \begin{bmatrix} 1 & 0 & 1 & 1 \\ 0 & 1 & 1 & 1 \\ 0 & 0 & \frac{\sqrt{3}}{\sqrt{2}} & \frac{\sqrt{3}}{\sqrt{2}} \\ 0 & 0 & -\frac{1}{\sqrt{2}} & \frac{1}{\sqrt{2}} \end{bmatrix}, \qquad J = \mathrm{diag}(1,1,-1,1), \qquad P = I.$$

If $R$ is partitioned in 4 (block) columns, the modulus strategy, in its first step, takes the first and the last column as the first pivot pair, and the middle two columns as the second pivot pair. The pivot matrices in these two processes are the same,

$$A_P^\circ = \begin{bmatrix} 1 & 1 \\ 1 & 1 \end{bmatrix},$$

and obviously singular. □

In the singular case, the Bunch-Parlett factorization gives the upper trapezoidal factor $R$ which is unsuitable for the trigonometric Jacobi method. The solution is either the two-sided Jacobi algorithm applied on $PH_P P^T$, or the QR factorization of $R^*$. In the latter case we obtain a full rank factor of a lower dimension. After the orthogonalization of this smaller factor, we have to assemble the full unitary matrix that diagonalizes $H_P$.

If the eigenvectors of $H$ are also needed, the algorithm requires the additional global matrix $U$, distributed in the same way as the block columns of $G^\circ = G^*$, starting with $U = I$. In each step, in addition to $G_P^\circ$, this matrix $U$ is locally multiplied by $U_P$ in order to accumulate the eigenvectors of $H$.

Instead of the Bunch-Parlett factorization in (2.2), the hyperbolic QR factorization (or JQR, for short) [18] can be applied, as well. Though a bit more accurate factor is produced, we are only interested in accumulating $U_P$ during the diagonalization, with columns as orthogonal as possible. Even with a factor of lower accuracy, $U_P$ can still be orthogonal to the machine precision, and, thus, useful for the approximate diagonalization of $H_P$. The JQR factorization, as already discussed, has a bookkeeping overhead, rendering it too inefficient, with no significant final accuracy gained.

The joint outline of all described parallel one-sided trigonometric algorithms is given as Algorithm 3.2.

For the pointwise algorithm, Dopico, Koev and Molera in [9, Theorem 7] showed that the computed eigenvector matrix $\widetilde{U}$ satisfies

$$\|\widetilde{U} - U\|_F = O(\epsilon \max\{n^{3/2} r, N_R\}),$$

where $n$ is the order of $H$, $r$ is the rank of $H$, while $N_R$ is the number of rotations used. If the algorithm is parallelized, under the standard IEEE model of the floating-point arithmetic, after $\ell$ stages, we obtain the linearized bound

$$\|\widetilde{U} - U\|_2 \leq 14\varepsilon\ell\sqrt{2pn}.$$

The proof of this fact can be found in the Appendix.

In principle, the eigenvector matrix $U$ can also be determined from the starting factor $G^*$, and the hyperbolic SVD of $G^*$, which can be written as

$$G^* U = V\Sigma. \tag{3.2}$$

Just note that $G^* U = V\Sigma$ is the final matrix computed by the algorithm. Multiplication of (3.2) from the left by $(V\Sigma)^* J$ yields

$$(V\Sigma)^* J G^* U = \Sigma^2 J = \Lambda, \tag{3.3}$$

where $\Lambda$ is a diagonal matrix of the computed eigenvalues. From (3.3) it follows that

$$U^* = \Lambda^{-1} (V\Sigma)^* J G^*.$$

**Algorithm 3.2:** Iterative part of trigonometric algorithms TB, TBC, TF, TFC

**Description**: Diagonalization of the matrix $H = GJG^*$ by the parallel one-sided trigonometric Jacobi algorithms. Assumption: $G^*$ is obtained by the Hermitian indefinite factorization of $H$ with complete pivoting, and then reordered to get $J$ partitioned as $J = \text{diag}(I_{npos}, -I_{n-npos})$. The matrix $G^*$ and the unitary matrix of eigenvectors $U$ are initially divided into $2p$ block columns, distributed so that the pair of blocks $(r+1, 2p-r)$ resides in process $r$. In each process the first block column is denoted by index $i$, and the second one by $j$.

```
Trigonometric_Jacobi(G, J, n);
begin
   repeat
        // compute the pivot submatrix H_P
      compute H_P from (3.1);
        // compute the Hermitian indefinite factorization of H_P
      Hermitian_Indefinite_Factorization_with_Complete_Pivoting(H_P,
      R);
      if algorithm = TB or TF then
         reorder the columns of R to their initial positions;
      end if
      if algorithm = block-oriented then
         if this is the first step in a sweep, annihilate all the off-diagonal elements
         of H_P;
         for all the other steps, annihilate only the elements of the block H_ij
         from (3.1);
         accumulate the matrix U_P of applied unitary transformations;
      else if algorithm = full block then
         diagonalize H_P;
         accumulate the unitary matrix U_P that diagonalizes H_p;
      end if
      if algorithm = TBC or TFC then
         apply the permutation from the Hermitian indefinite factorization to the
         rows of U_P;
      end if
      [G_i^* G_j^*] = [G_i^* G_j^*] · U_P;
        // accumulate eigenvectors
      [U_i U_j] = [U_i U_j] · U_P;
      send/receive one of the blocks in [G_i^* G_j^*] and [U_i U_j] to/from the
      neighboring process according to the modulus strategy;
   until convergence;
end
```

However, the relative accuracy of such a computation of eigenvectors (i.e., a bound on angles between the computed eigenvectors and the corresponding eigenvector subspaces) has not yet been proven, nor extensively tested.

## 3.2 Hyperbolic Jacobi algorithms

Hyperbolic Jacobi algorithms, instead of only trigonometric rotations, use both trigonometric and hyperbolic rotations, depending on the signs in $J$. The algorithms diagonalize, by congruence, a definite matrix pair $(A^\circ, J) := (G^*G, J)$.

Each process owns two block columns, $G_i$ and $G_j$, and the pivot block is

$$A_P^\circ := A_P = \begin{bmatrix} A_{ii} & A_{ij} \\ A_{ij}^* & A_{jj} \end{bmatrix} = \begin{bmatrix} G_i^* \\ G_j^* \end{bmatrix} \begin{bmatrix} G_i & G_j \end{bmatrix}. \tag{3.4}$$

Let $J_P$ be a part of $J$ corresponding to $A_P^\circ$. Then the definite matrix pair $(A_P^\circ, J_P)$ is transformed per process.

Due to the full column rank of the initial matrix $G^\circ = G$, the matrix $[G_i, G_j]$ also has full column rank, and $A_P$ is a positive definite matrix. Thus, it can be factored by the Cholesky factorization (with or without pivoting):

$$P^T A_P P = R^* R.$$

The rest of the computation is very similar to the trigonometric case. We need to apply the sequence of rotations to the columns of $R$. When the signs of diagonal elements in $J_P$ match, a trigonometric rotation is applied, otherwise a hyperbolic one is used. The parameters of a hyperbolic rotation are computed as described in [22]. Despite being impossible in theory, $\tanh \varphi = \pm 1$ can occur in finite precision arithmetic. In these extremely rare cases, a smaller angle $\varphi'$ is used for the hyperbolic transformation (usually, $\tanh \varphi' = \pm 0.9$).

Here, the rotations are accumulated, per step, in the local matrix $V_P^{-*}$, to be applied later to the block $G_P^\circ$, but the global $V^{-*}$ no longer needs to be maintained, since $GV^{-*} = U\Sigma$. The final eigenvector matrix $U$ is obtained by normalizing the columns of the resulting matrix $GV^{-*}$.

The hyperbolic algorithms explained so far, the full block (HF) and the block-oriented (HB), can be tuned even further. If we follow the idea from the trigonometric case (presented in [7] for the Jacobi SVD algorithm) and use the Cholesky factorization with diagonal pivoting, some speedup is gained in the full block case (HFC), but lost it in the block-oriented case (HBC). Similarly to the sorted trigonometric case (TBC, TFC), the eigenvalues are computed in decreasing order by their absolute value. This happens because "sorting" (by pivoting) can spoil the quadratic convergence by mixing columns with almost the same hyperbolic norms, which correspond to different signs in $J$ (see [19] for details).

On the other hand, we can do the Cholesky factorization in two parts, respecting the positive and negative signs in $J$. Suppose that $A_P$ has the following block structure

$$A_P = \begin{bmatrix} A_{11} & A_{12} \\ A_{12}^* & A_{22} \end{bmatrix}, \tag{3.5}$$

and the square diagonal blocks $A_{11}$ and $A_{22}$ correspond to the positive/negative signs in $J_P$. If $A_{22}$ in (3.5) does not exist ($J_P = I$), the whole block $A_P$ is factored by the Cholesky factorization with diagonal pivoting. If $A_{11}$ is non-existent ($J_P = -I$), the same is done, but afterwards the columns are reversed to keep the column norms in non-decreasing order. Else, the block $A_{11}$ is factored by the Cholesky factorization with diagonal pivoting, $P_1^T A_{11} P_1 = R_{11}^* R_{11}$, so (3.5) becomes

$$\begin{bmatrix} P_1^T & \\ & I \end{bmatrix} \begin{bmatrix} A_{11} & A_{12} \\ A_{12}^* & A_{22} \end{bmatrix} \begin{bmatrix} P_1 & \\ & I \end{bmatrix} = \begin{bmatrix} R_{11}^* & \\ R_{12}^* & I \end{bmatrix} \begin{bmatrix} I & \\ & S \end{bmatrix} \begin{bmatrix} R_{11} & R_{12} \\ & I \end{bmatrix}. \tag{3.6}$$

From (3.6) it follows that $R_{11}^* R_{12} = P_1^T A_{12}$, and $R_{12}$ can be computed by solving this triangular linear system. The Schur complement
$$S = A_{22} - R_{12}^* R_{12}$$
is also a Hermitian positive definite matrix, so we can factorize it by the Cholesky factorization with diagonal pivoting
$$P_2^T S P_2 = R_{22}^* R_{22}. \tag{3.7}$$
Combining (3.7) and (3.6), we have
$$\begin{bmatrix} P_1^T & \\ & I \end{bmatrix} A_P \begin{bmatrix} P_1 & \\ & I \end{bmatrix} = \begin{bmatrix} R_{11}^* & \\ R_{12}^* & I \end{bmatrix} \begin{bmatrix} I & \\ & P_2 \end{bmatrix} \begin{bmatrix} I & \\ & P_2^T S P_2 \end{bmatrix} \begin{bmatrix} I & \\ & P_2^T \end{bmatrix} \begin{bmatrix} R_{11} & R_{12} \\ & I \end{bmatrix}$$
$$= \begin{bmatrix} I & \\ & P_2 \end{bmatrix} \begin{bmatrix} R_{11}^* & \\ P_2^T R_{12}^* & R_{22}^* \end{bmatrix} \begin{bmatrix} R_{11} & R_{12} P_2 \\ & R_{22} \end{bmatrix} \begin{bmatrix} I & \\ & P_2^T \end{bmatrix}.$$
The final sign-pivoted Cholesky factorization of $A_P$ is given by
$$\begin{bmatrix} P_1^T & \\ & P_2^T \end{bmatrix} A_P \begin{bmatrix} P_1 & \\ & P_2 \end{bmatrix} = R_S^* R_S, \qquad R_S := \begin{bmatrix} R_{11} & R_{12} P_2 \\ & R_{22} \end{bmatrix}.$$
The block column structure of the factor $R_S$ can be written as
$$R_S = \begin{bmatrix} R_{S,1} & R_{S,2} \end{bmatrix}, \qquad R_{S,1} = \begin{bmatrix} R_{11} \\ 0 \end{bmatrix}, \qquad R_{S,2} = \begin{bmatrix} R_{12} P_2 \\ R_{22} \end{bmatrix}.$$

After the factorization, we need to reverse the order of columns in $R_{S,2}$ — the first column is swapped with the last, the second one with the penultimate, and so on. This reversal tries to maintain the eigenvalues sorted in non-decreasing order, thus, ensuring the quadratic convergence in the case of multiple or clustered eigenvalues [11]. Moreover, the computed eigenvalues are sorted "for free".

Algorithm 3.3 computes the specially pivoted Cholesky factorization used in algorithms HBSC and HFSC.

The joint outline of all described parallel one-sided hyperbolic Jacobi algorithms is given in Algorithm 3.4. The full block (HFSC) and the block-oriented (HBSC) algorithms, are faster than their respective diagonally pivoted counterparts HFC and HBC, the former one more than the latter.

## 4 Numerical testing

Numerical testing has been conducted on a cluster of 16 blade servers, each equipped with two dual-core Intel Xeon 5150 CPUs at 2.66 GHz, and with 8 GB of RAM.

The software used consists of Intel Fortran and C++ compilers and Math Kernel Library 11.0.074 for EM64T on GNU/Linux, and Open MPI 1.3. Our programs are single-threaded and, mostly, Fortran 77 and 90 based.

We have tested nonsingular real symmetric matrices of orders $n$ from 1000 to 16000, in steps of 1000. In each test matrix, the elements of the upper triangle have been pseudorandomly generated, with the uniform distribution of elements in $[-5, 5]$, by using the LAPACK testing routine `dlarnd`. This kind of generation produces matrices with tightly clustered eigenvalues, varying from $10^{-3}$ to $10^3$ in magnitude.

Finally, the input matrices $G$ and $J$ have been computed by the sequential Bunch-Parlett factorization with complete pivoting. The generation and preprocessing of matrices is not included in the times given below.

A performance comparison of all described trigonometric and hyperbolic parallel Jacobi algorithms is given in Fig. 4 and Fig. 5, respectively. The whole computation has been performed in IEEE double precision arithmetic.

At the first glance more natural, the trigonometric Jacobi algorithms are, in fact, slower and less accurate than the hyperbolic ones. For large matrices, the orthogonality of the

**Algorithm 3.3:** Specially pivoted Cholesky factorization with respect to $J_P$

**Description**: The algorithm computes the specially pivoted Cholesky factorization of the block matrix $A_P$ of order $n_P$, according to the number of positive signs *npos* in $J_P$. The computed factor $R_S$ is returned in $A_P$, while the permutation $P$ is returned in a separate vector.

**Input**: $n_P$, *npos*.
**Input/Output**: $A_P$.
**Output**: $P$.

```
J_P_Pivoted_Cholesky(A_P, n_P, npos, P);
```
**begin**
    **if** *npos* $= n_P$ **then**    // case $J = I$
        factorize $A_{11}$ by the Cholesky factorization with diagonal pivoting;
        **return** $R$, $P$
    **else if** *npos* $= 0$ **then**    // case $J = -I$
        factorize $A_{22}$ by the Cholesky factorization with diagonal pivoting;
        reverse the columns of $R$ and $P$;
        **return** $R$, $P$
    **else**    // case $J \neq \pm I$
        factorize $A_{11}$ by the Cholesky factorization with diagonal pivoting;
        compute $R_{12}$ from the triangular linear system $R_{11}^* R_{12} = P_1^T A_{12}$;
        $S = A_{22} - R_{12}^* R_{12}$;
        factorize $S$ by the Cholesky factorization with diagonal pivoting,
        $P_2^T S P_2 = R_{22}^* R_{22}$;
        apply the permutation $P_2$ from the right to $R_{12}$;
        reverse the columns of $R_{S,2}$ and $P_2$;
        **return** $R_S$, $P = \mathrm{diag}(P_1, P_2)$
    **end if**
**end**

**Algorithm 3.4:** Iterative part of hyperbolic algorithms HB, HBC, HBSC, HF, HFC, HFSC

**Description**: Diagonalization of the pair $(A, J) = (G^*G, J)$ by the parallel one-sided hyperbolic Jacobi algorithms. Assumption: $G$ is obtained by the Hermitian indefinite factorization of $H$ with complete pivoting, and then reordered to get $J$ partitioned as $J = \text{diag}(I_{npos}, -I_{n-npos})$. The matrix $G$ is initially divided into $2p$ block columns, distributed so that the pair of blocks $(r+1, 2p-r)$ resides in process $r$. In each process the first block column is denoted by index $i$, and the second one by $j$.

```
Hyperbolic_Jacobi(G, J, n);
begin
    repeat
        // compute the pivot submatrix A_P
        compute A_P from (3.4);
        J_P = diag(J_ii, J_jj);
        // compute a selected type of the Cholesky factorization of
            A_P
        if algorithm = HBSC or HFSC then
            J_P_Pivoted_Cholesky(A_P, R);
        else if algorithm = HBC or HFC then
            Diagonally_Pivoted_Cholesky_with_Reordering(A_P, R);
        else
            Unpivoted_Cholesky(A_P, R);
        end if
        if algorithm = block-oriented then
            if this is the first step in a sweep, annihilate all the off-diagonal elements
            of A_P;
            for all the other steps, annihilate only the elements of the block A_ij
            from (3.4);
            accumulate the J_P-unitary matrix V_P^{-*} of applied transformations;
        else if algorithm = full block then
            diagonalize A_P;
            accumulate the J_P-unitary matrix V_P^{-*} that diagonalizes the pair
            (R^*R, J_P);
        end if
        if algorithm = HBC, HFC, HBSC or HFSC then
            apply the permutation from the Cholesky factorization to the rows of
                V_P^{-*};
        end if
        [G_i G_j] = [G_i G_j] · V_P^{-*};
        send/receive one of the blocks in [G_i G_j] to/from the neighboring process
        according to the modulus strategy;
    until convergence;
end
```

(a) full block

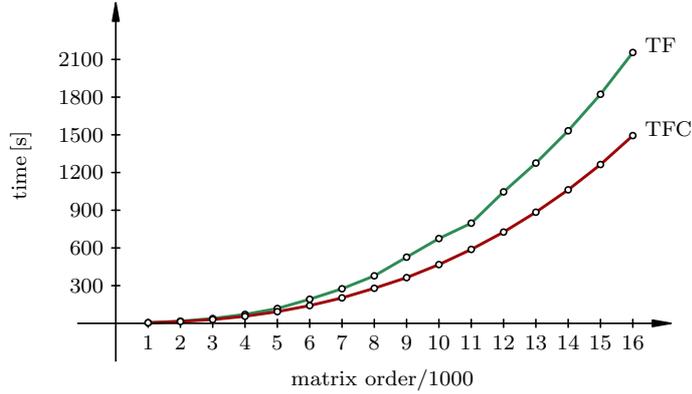

(b) block-oriented

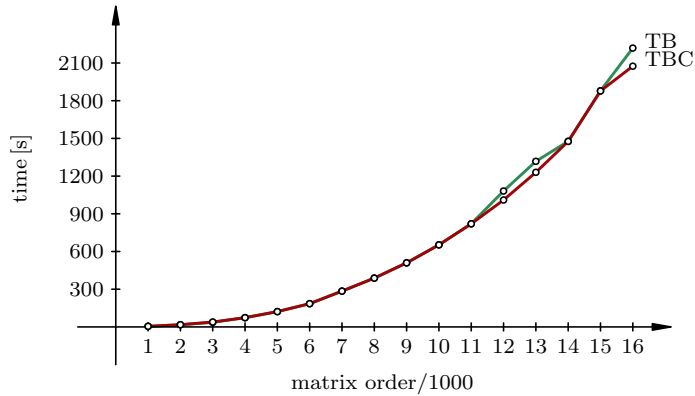

Figure 4: Comparison between trigonometric Jacobi algorithms on 64 processes.

computed eigenvectors $\widetilde{U}$, measured by $\|\widetilde{U}\widetilde{U}^T - I\|_F$ and $\|\widetilde{U}^T\widetilde{U} - I\|_F$, is approximately $10^{-15}$ or less in the hyperbolic case, while it is not below $10^{-12}$ in the trigonometric case. We believe that the repeated accumulation of $U$ in each step contributes significantly to the overall error, while slowing down the algorithm.

Our conclusion is that the hyperbolic parallel Jacobi algorithms are *easier* to implement, and could be substantially *faster* than the trigonometric ones.

The major speedup — dramatically reducing the number of sweeps, is obtained by the special pivoting, both in the Bunch-Parlett (for trigonometric) and in the Cholesky factorization (for hyperbolic algorithms). This effect is more visible in the full block case, and somewhat less in the block-oriented case. For example, on matrices of order 16000, the number of in-process sweeps is reduced from 17 in the "classical" HF implementation, to 11 in the HFSC algorithm, resulting in 30% speedup.

Therefore, we expect that, given a similar computing architecture and matrices large enough, the performance of HFSC and HBSC should be superior to the other algorithms described herein.

The timings shown in Fig. 4 and Fig. 5 are accurate in the sense that the successive runs of the same algorithm on the same data differ by less than a second, mostly due to the fact that all MPI communication is synchronous in our implementations. Moreover, for performance reasons, the threads should be assigned to only one processor core during the entire run. Otherwise, a thread could be rescheduled at any time to an another core, at the operating system's discretion. This move causes a severe cache invalidation and slowdown, which propagates through the entire run.

Fig. 6 shows the portion of the total execution time spent in communication. The

(a) full block

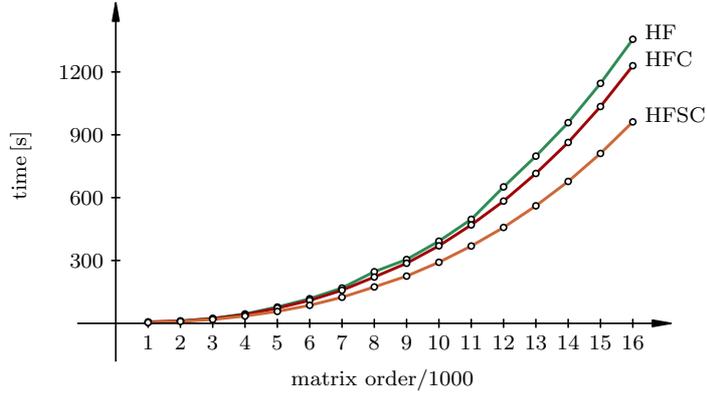

(b) block-oriented

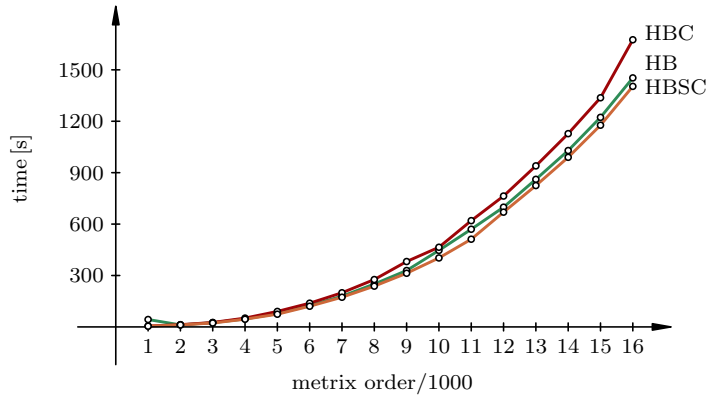

Figure 5: Comparison between hyperbolic Jacobi algorithms on 64 processes.

relatively high percentage is caused by the synchronous nature of the communication and synchronization — the observed execution time of a step is always the time of an MPI process with the slowest computation in that step.

This also exposes a side-effect of dealing with indefinite matrices. If one block is definite, and the other one indefinite, their diagonalization times may differ by up to 30%, which hurts the full block algorithms. It remains an open question whether the asynchronous communication patterns could hide these latencies in part.

We have also compared our algorithms with the ScaLAPACK routine `pdsygvx`, which has been used to compute the eigenvalues of the real definite matrix pair $(A, J)$, where $A := G^*G$ is a positive definite matrix (see Fig. 7). To get a fair comparison, the computation of $A$ from $G$ (`pdsyrk`) is *not* added to the total time.

This comparison has been done on a cluster similar to the described one, with the main difference being the 10 Gb Ethernet interconnection. The parameters of `pdsygvx` have been chosen as recommended ones for the most accurate eigenvalue computation. We have experimentally found that $64 \times 64$ blocks is the fastest blocking option.

The one-sided Jacobi algorithms inside each process can be replaced by other, faster algorithms, provided that these algorithms give numerically ($J$–)orthogonal eigenvectors in the trigonometric (hyperbolic) case. Since the outer Jacobi method is self-correcting, these eigenvectors need not be very close to the exact ones.

For example, in the full block trigonometric case, the inner (per block) Jacobi algorithm can be replaced by the tridiagonalization and divide-and-conquer algorithm. Such a modification (TFDC, for short) has been tested with the LAPACK `dsyevd` routine. The total running times for TFC and TFDC are given in Table 1.

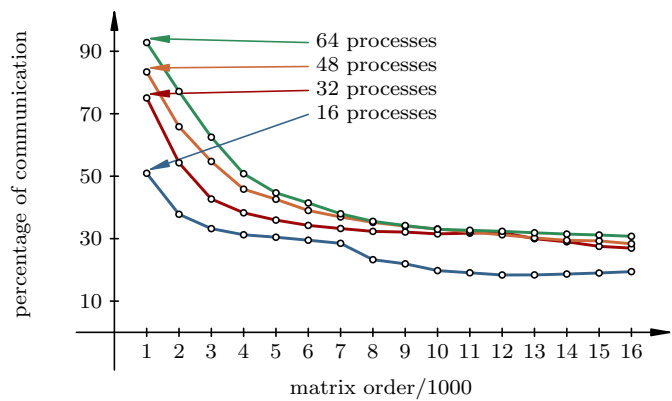

(a) full block (HFSC)

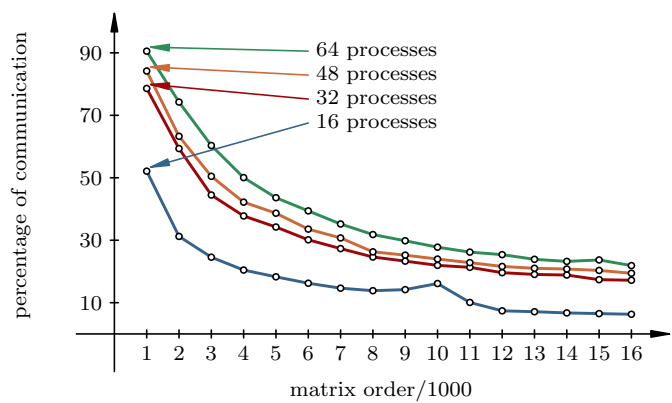

(b) block-oriented (HBSC)

Figure 6: Percentage of total time spent in communication.

While the stopping criterion of the one-sided Jacobi algorithm is natural, there is no easy way to stop the whole algorithm when per block divide-and-conquer algorithm is employed. To obtain the comparison in Table 1, if the one-sided trigonometric full block algorithm uses $s$ sweeps to convergence, we run $s - 2$ sweeps with per block divide-and-conquer algorithm, and the final sweeps (until convergence) with the Jacobi algorithm. A more detailed study is needed to determine the optimal stopping criterion with per block divide-and-conquer algorithm.

# Conclusion

The presented parallel algorithms can be even faster, if we apply the sequential blocking for large matrices inside each process. These algorithms, called the three-level Jacobi algorithms, are the subject of another paper [20]. The switching point between the non-blocked and the blocked algorithm inside each process depends on various factors, such as the processor speed and organization, and the speed of the interconnection network. Once the hardware is fixed, profiling (like in the self-tuning packages, e.g., ATLAS) is needed to reveal that switching point for each algorithm.

# A  Appendix

Here we present a floating-point error analysis for the accumulation of products of nearly orthogonal matrices, which is used to compute the eigenvectors in the trigonometric Jacobi

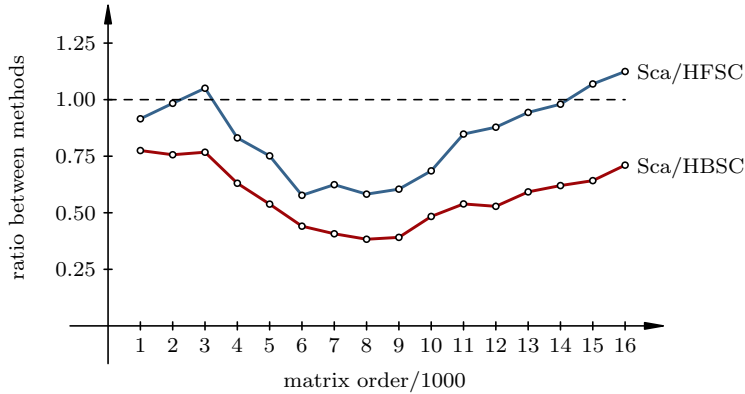

Figure 7: Time ratio between ScaLAPACK and HFSC, HBSC on 64 processes.

|  |  | Matrix sizes | | | |
| --- | --- | --- | --- | --- | --- |
| Routine | processes | 4000 | 8000 | 12000 | 16000 |
| TFC | 32 | 67 | 369 | 1034 | 2525 |
| TFDC | 32 | 60 | 316 | 849 | 1796 |
| TFC | 64 | 53 | 298 | 749 | 1492 |
| TFDC | 64 | 45 | 254 | 697 | 1373 |

Table 1: Timing (in seconds) of TFC and TFDC algorithms.

algorithms.

We use the IEEE standard model of floating-point arithmetic

$$\mathrm{fl}(a \circ b) = (a \circ b)(1 + \varepsilon_\circ), \quad |\varepsilon_\circ| \leq \varepsilon,$$

where $\circ$ is any of the four arithmetic operations, and $\varepsilon$ is the unit roundoff error. Additionally, we also assume that the square roots can be computed with the same accuracy. In our analysis, numerically subscribed $\varepsilon$'s denote the quantities bounded by the unit roundoff error.

In a pivot block which is assigned to a particular process, in each inner step, we choose and orthogonalize a pair of columns. For simplicity, we may assume that these inner steps occur simultaneously in all processes. This single simultaneous transformation of $p$ pairs of columns will be called a *stage*.

At the beginning of the algorithm, the eigenvector matrix is set to the identity matrix, i.e., $U^{(0)} := I$. The superscript denotes the stage index throughout the iterations, and the transformation $U^{(\ell-1)} \mapsto U^{(\ell)}$ in stage $\ell$ consists of $p$ independent trigonometric rotations. The $i$-th and $j$-th columns of the matrix $U^{(\ell)}$ are denoted by $u_i^{(\ell)}$ and $u_j^{(\ell)}$, respectively. In addition, all the computed quantities are denoted by tilde, i.e., $\tilde{u}_i^{(\ell)}$ and $\tilde{u}_j^{(\ell)}$ denote the $i$-th and $j$-th columns of the computed matrix $\widetilde{U}^{(\ell)}$.

Let $i$ and $j$ be the indices of columns transformed in one of the processes in stage $\ell$. When the columns are transformed by a trigonometric rotation, the exactly computed new columns would be

$$\left[u_i^{(\ell)}, u_j^{(\ell)}\right] = \left[u_i^{(\ell-1)}, u_j^{(\ell-1)}\right] \begin{bmatrix} c & s \\ -s & c \end{bmatrix} = \left[cu_i^{(\ell-1)} - su_j^{(\ell-1)}, su_i^{(\ell-1)} + cu_j^{(\ell-1)}\right]. \quad \text{(A.1)}$$

Instead of (A.1), the new columns actually computed in the floating-point arithmetic are

$$\left[\tilde{u}_i^{(\ell)}, \tilde{u}_j^{(\ell)}\right] = \left[\tilde{u}_i^{(\ell-1)}, \tilde{u}_j^{(\ell-1)}\right] \begin{bmatrix} \tilde{c} & \tilde{s} \\ -\tilde{s} & \tilde{c} \end{bmatrix}. \quad \text{(A.2)}$$

Here, we assume that the columns $\tilde{u}_i^{(\ell-1)}$ and $\tilde{u}_j^{(\ell-1)}$ have already accumulated some error from the previous transformations. In the formula (A.2), $\tilde{c}$ and $\tilde{s}$ denote the computed cosine and sine, respectively, so

$$\tilde{c} = (1 + \varepsilon_c)c, \quad \tilde{s} = (1 + \varepsilon_s)s.$$

Typically (see, for example, [9, equation (30)]), we have

$$|\varepsilon_c| \leq \frac{5\varepsilon}{1 - 5\varepsilon}, \quad |\varepsilon_s| \leq \frac{5\varepsilon}{1 - 5\varepsilon}. \tag{A.3}$$

The computed elements of the columns $\tilde{u}_i^{(\ell)}$ and $\tilde{u}_j^{(\ell)}$ at the end of stage $\ell$ are

$$\tilde{u}_{ki}^{(\ell)} = \Big[(1+\varepsilon_c)(1+\varepsilon_1)c\tilde{u}_{ki}^{(\ell-1)} - (1+\varepsilon_s)(1+\varepsilon_2)s\tilde{u}_{kj}^{(\ell-1)}\Big](1+\varepsilon_3)$$

$$\tilde{u}_{kj}^{(\ell)} = \Big[(1+\varepsilon_s)(1+\varepsilon_4)s\tilde{u}_{ki}^{(\ell-1)} + (1+\varepsilon_c)(1+\varepsilon_5)c\tilde{u}_{kj}^{(\ell-1)}\Big](1+\varepsilon_6).$$

In other words,

$$\tilde{u}_i^{(\ell)} = \big(c\tilde{u}_i^{(\ell-1)} - s\tilde{u}_j^{(\ell-1)}\big) + \delta u_i^{(\ell)}, \quad \tilde{u}_j^{(\ell)} = \big(s\tilde{u}_i^{(\ell-1)} + c\tilde{u}_j^{(\ell-1)}\big) + \delta u_j^{(\ell)},$$

and the elementwise errors, linearized up to the term of order $O(\varepsilon^2)$, are

$$\begin{aligned}
\delta u_{ki}^{(\ell)} &= \varepsilon^{(1)} c\tilde{u}_{ki}^{(\ell-1)} - \varepsilon^{(2)} s\tilde{u}_{kj}^{(\ell-1)}, & \varepsilon^{(1)} &= \varepsilon_c + \varepsilon_1 + \varepsilon_3, & \varepsilon^{(2)} &= \varepsilon_s + \varepsilon_2 + \varepsilon_3, \\
\delta u_{kj}^{(\ell)} &= \varepsilon^{(3)} s\tilde{u}_{ki}^{(\ell-1)} + \varepsilon^{(4)} c\tilde{u}_{kj}^{(\ell-1)}, & \varepsilon^{(3)} &= \varepsilon_s + \varepsilon_4 + \varepsilon_6, & \varepsilon^{(4)} &= \varepsilon_c + \varepsilon_5 + \varepsilon_6,
\end{aligned} \tag{A.4}$$

for $k = 1, \ldots, n$.

The previous formula can be written in a matrix form, which summarizes all $p$ transformations in stage $\ell$. Instead of the orthogonal matrix $U^{(\ell)}$, we have computed a slightly perturbed matrix $\widetilde{U}^{(\ell)}$

$$\widetilde{U}^{(\ell)} = \mathrm{fl}(\widetilde{U}^{(\ell-1)}\widetilde{Q}^{(\ell)}) = \widetilde{U}^{(\ell-1)}Q^{(\ell)} + \delta U^{(\ell)}, \tag{A.5}$$

where $\widetilde{Q}^{(\ell)}$ is the computed matrix of $p$ independent rotations used in stage $\ell$, while $Q^{(\ell)}$ is its exact and orthogonal counterpart. So, (A.5) is a matrix equivalent of (A.2) for all processes.

**Theorem 1.** *Let $U^{(\ell)}$ be the exact orthogonal matrix of accumulated transformations after $\ell$ stages of the parallel trigonometric Jacobi algorithm, and let $\widetilde{U}^{(\ell)}$ be the computed matrix in the floating-point arithmetic. Then*

$$\|\widetilde{U}^{(\ell)} - U^{(\ell)}\|_2 \leq \sum_{m=1}^{\ell} \|\delta U^{(m)}\|_2, \tag{A.6}$$

*where $\delta U^{(m)}$ is the perturbation matrix defined by (A.5) in stage $m$ of the algorithm. Moreover, if the error in all computed cosines and sines is bounded by (A.3), then*

$$\|\widetilde{U}^{(\ell)} - U^{(\ell)}\|_2 \leq 14\varepsilon\ell\sqrt{2pn}, \tag{A.7}$$

*with the right-hand side linearized up to the term of order $O(\varepsilon^2)$.*

*Proof.* The proof follows by induction over the stage index $\ell$. Let $Z^{(\ell)}$ be the error committed after $\ell$ stages of the algorithm,

$$Z^{(\ell)} := \widetilde{U}^{(\ell)} - U^{(\ell)},$$

for any $\ell \geq 0$. At the beginning of the first stage, we have $\widetilde{U}^{(0)} = U^{(0)} = I$, and $Z^{(0)} = 0$, so both claims are obviously true for $\ell = 0$.

Now suppose that the claim is valid at the beginning of stage $\ell \geq 1$, i.e., the computed matrix accumulated so far is $\widetilde{U}^{(\ell-1)} = U^{(\ell-1)} + Z^{(\ell-1)}$, with

$$\|Z^{(\ell-1)}\|_2 \leq \sum_{m=1}^{\ell-1} \|\delta U^{(m)}\|_2. \tag{A.8}$$

After the transformation in stage $\ell$ (in all processes), the computed new matrix is $\widetilde{U}^{(\ell)} = U^{(\ell)} + Z^{(\ell)}$, and from (A.5) it follows that

$$\begin{aligned}
\widetilde{U}^{(\ell)} &= \widetilde{U}^{(\ell-1)} Q^{(\ell)} + \delta U^{(\ell)} \\
&= \bigl(U^{(\ell-1)} + Z^{(\ell-1)}\bigr) Q^{(\ell)} + \delta U^{(\ell)} \\
&= U^{(\ell-1)} Q^{(\ell)} + Z^{(\ell-1)} Q^{(\ell)} + \delta U^{(\ell)}.
\end{aligned}$$

Since $U^{(\ell)} = U^{(\ell-1)} Q^{(\ell)}$ is the exact new matrix after $\ell$ stages (a matrix equivalent of (A.1)), the new error matrix $Z^{(\ell)}$ can be written as

$$Z^{(\ell)} = Z^{(\ell-1)} Q^{(\ell)} + \delta U^{(\ell)}.$$

The exact transformation matrix $Q^{(\ell)}$ is orthogonal, and by unitary invariance of the spectral norm, we get

$$\|Z^{(\ell)}\|_2 \leq \|Z^{(\ell-1)}\|_2 + \|\delta U^{(\ell)}\|_2, \tag{A.9}$$

By the induction hypothesis (A.8), this completes the proof of (A.6).

To prove the second claim, we need to bound the perturbation $\delta U^{(\ell)}$ in (A.5). We shall prove a slightly more general result, and (A.7) will then follow by using the bounds in (A.3).

We assume that all the cosines and sines in the algorithm are computed with relative errors $\varepsilon_c$ and $\varepsilon_s$ bounded by

$$|\varepsilon_c| \leq \varepsilon_{\cos}, \quad |\varepsilon_s| \leq \varepsilon_{\sin}, \tag{A.10}$$

where $\varepsilon_{\cos}$ and $\varepsilon_{\sin}$ depend on a particular algorithm that is used to compute the elements of each trigonometric rotation. The only requirement is that both bounds must satisfy

$$\varepsilon_{\cos} = O(\varepsilon), \quad \varepsilon_{\sin} = O(\varepsilon), \tag{A.11}$$

linearized up to the terms of order $O(\varepsilon^2)$, with small "hidden" constants, which may be different. In this setting, the bound in (A.7) becomes

$$\|\widetilde{U}^{(\ell)} - U^{(\ell)}\|_2 \leq (\varepsilon_{\cos} + \varepsilon_{\sin} + 4\varepsilon)\ell\sqrt{2pn}, \tag{A.12}$$

which is again true for $\ell = 0$.

Suppose that, at the beginning of stage $\ell \geq 1$, the error $Z^{(\ell-1)}$ in the computed matrix $\widetilde{U}^{(\ell-1)}$ is bounded by (A.12), i.e.,

$$\|Z^{(\ell-1)}\|_2 \leq (\varepsilon_{\cos} + \varepsilon_{\sin} + 4\varepsilon)(\ell-1)\sqrt{2pn}. \tag{A.13}$$

From $\widetilde{U}^{(\ell-1)} = U^{(\ell-1)} + Z^{(\ell-1)}$ it immediately follows that

$$\|\widetilde{U}^{(\ell-1)}\|_2 \leq \|U^{(\ell-1)}\|_2 + \|Z^{(\ell-1)}\|_2 = 1 + \|Z^{(\ell-1)}\|_2.$$

Thus, by (A.11) and (A.13), all the elements of the perturbed matrix $\widetilde{U}^{(\ell-1)}$ satisfy

$$|\tilde{u}_{kl}^{(\ell-1)}| \leq 1 + O(\varepsilon), \quad k, l = 1, \ldots, n.$$

Since $|c|, |s| \leq 1$, from (A.4) it follows that the elementwise perturbations in the transformed columns $i$ and $j$ (in a particular process) can be bounded by

$$|\delta u_{ki}^{(\ell)}|, |\delta u_{kj}^{(\ell)}| \leq \bigl(|\varepsilon_c| + |\varepsilon_s| + 4\varepsilon\bigr)\bigl(1 + O(\varepsilon)\bigr).$$

Of course, if a column is not transformed in stage $\ell$, the corresponding perturbations are equal to zero. Let $E^{(\ell)}$ be a matrix of order $n$ with the following column structure — a column of $E^{(\ell)}$ is equal to zero if this column is not transformed in stage $\ell$, otherwise it is equal to $e$, where $e$ is a vector with all elements equal to one. Since exactly $2p$ columns are transformed in stage $\ell$, the matrix $E^{(\ell)}$ has $2p$ columns equal to $e$.

Let the symbol $|\ |$ denote the pointwise absolute value of a matrix. From the above argument, by using (A.10) over all $p$ processes, we get the following linearized bound for the perturbation $\delta U^{(\ell)}$ in (A.5)

$$|\delta U^{(\ell)}| \leq (\varepsilon_{\cos} + \varepsilon_{\sin} + 4\varepsilon)E^{(\ell)}.$$

For any two matrices $A, B \in \mathbb{R}^{n \times n}$, if $|A| \leq B$, then $\|A\|_2 \leq \|B\|_2$ (see [15, Lemma 6.6.(b)]). By taking $A = \delta U^{(\ell)}$, and $B = (\varepsilon_{\cos} + \varepsilon_{\sin} + 4\varepsilon)E^{(\ell)}$, it follows that

$$\|\delta U^{(\ell)}\|_2 \leq (\varepsilon_{\cos} + \varepsilon_{\sin} + 4\varepsilon)\|E^{(\ell)}\|_2.$$

Since $E^{(\ell)}$ is of rank one, we have $\|E^{(\ell)}\|_2^2 = \text{trace}([E^{(\ell)}]^T E^{(\ell)}) = 2pn$, so

$$\|\delta U^{(\ell)}\|_2 \leq (\varepsilon_{\cos} + \varepsilon_{\sin} + 4\varepsilon)\sqrt{2pn}.$$

Now, (A.12) follows immediately from (A.9) and (A.13).

Finally, if we take the linearized bounds from (A.3), $\varepsilon_{\cos} = \varepsilon_{\sin} = 5\varepsilon$, the relation (A.12) becomes

$$\|\widetilde{U}^{(\ell)} - U^{(\ell)}\|_2 \leq 14\varepsilon\ell\sqrt{2pn},$$

which proves (A.7). □ □

In the sequential case $p = 1$, the bound (A.7) in the spectral norm is similar to the perturbation bound (in the Frobenius norm) given by Dopico, Koev and Molera in [9, Theorem 6].

The departure from orthonormality of the computed eigenvector matrix $\widetilde{U}$ in the spectral norm can also be estimated by using the bound (A.7). We have

$$\|\widetilde{U}^T\widetilde{U} - I\|_2 = \|\widetilde{U}^T\widetilde{U} - U^T U\|_2 = \|\widetilde{U}^T\widetilde{U} - U^T\widetilde{U} + U^T\widetilde{U} - U^T U\|_2$$
$$= \|(\widetilde{U} - U)^T\widetilde{U} + U^T(\widetilde{U} - U)\|_2 \leq \|\widetilde{U} - U\|_2 \cdot \|\widetilde{U}\|_2 + 1 \cdot \|\widetilde{U} - U\|_2. \quad \text{(A.14)}$$

Since

$$\|\widetilde{U}\|_2 = \|\widetilde{U} - U + U\|_2 \leq \|\widetilde{U} - U\|_2 + \|U\|_2 = \|\widetilde{U} - U\|_2 + 1,$$

then (A.14) simplifies to

$$\|\widetilde{U}^T\widetilde{U} - I\|_2 \leq \|\widetilde{U} - U\|_2(\|\widetilde{U} - U\|_2 + 2).$$

Now, if the total number of stages is known, one can use this and the formula (A.7) to bound the departure from orthonormality.

## Acknowledgements

We wish to thank Prof. Dr. Mioara Mandea, Dr. Vincent Lesur, and Alexander Jordan from Department of Earth's Magnetic Field, Helmholtz Centre, Potsdam, for providing us with testing time on their cluster. We also thank Dr. Vedran Šego for his valuable suggestions and proof-reading.